\documentclass[a4paper,10pt]{article}
\usepackage{amsmath,amsthm,amssymb}
\usepackage{graphicx,subfigure}
\usepackage{url}

\usepackage{algorithm2e}
\usepackage{listings}
\usepackage{color}

\usepackage{multicol}
\def\twoplot[#1]#2#3#4#5{
\begin{figure}[h]
\begin{multicols}{2}
\begin{center}
    \includegraphics*[#1]{#2}
    \caption{\label{#2} #4}
\end{center}
\begin{center}
    \includegraphics*[#1]{#3}
    \caption{\label{#3} #5}
\end{center}
\end{multicols}
\end{figure}
}% end twoplot macro

\usepackage[colorlinks,bookmarksopen,bookmarksnumbered,citecolor=red,urlcolor=red]{hyperref}
\hypersetup{pdftitle={A shallow water with variable pressure model for blood flow simulation},
bookmarks=true,
pdftoolbar=true,
pdfmenubar=true,
pdfauthor={O. Delestre, A. R. Ghigo, J.-M. Fullana, P.-Y. Lagr\'ee},
pdfsubject={blood flow, finite volume, numerical method, well balanced scheme, HLL flux, shallow water equations},
pdfcreator={Delestre},
pdfproducer={Delestre},
pdfkeywords={blood flow}{shallow water}{finite volume}{well-balanced}{HLL flux}{hydrostatic reconstruction}{source term}
{steady state}{numerical method}}

\title{A shallow water with variable pressure model for blood flow simulation}

\author{{O. Delestre}\footnote{Lab. J.A. Dieudonn\'e \& EPU Nice Sophia, University of
 Nice, France, e-mail : delestre@math.unice.fr}, A.R. Ghigo\footnote{CNRS and UPMC Universit\'e Paris 06, UMR 7190,
 Institut Jean Le Rond d'Alembert, France, e-mail : arthur.ghigo@dalembert.upmc.fr}, {J.-M. Fullana}\footnote{CNRS and UPMC
 Universit\'e Paris 06, UMR 7190, Institut Jean Le Rond d'Alembert, France, e-mail : fullana@lmm.jussieu.fr}
and {P.-Y. Lagr\'ee}\footnote{CNRS and UPMC Universit\'e Paris 06, UMR 7190, Institut Jean Le Rond d'Alembert, France, e-mail :
 pierre-yves.lagree@upmc.fr}}

\begin{document}
\maketitle

\begin{abstract} We performed numerical simulations of blood flow in arteries with a variable stiffness and cross-section at rest
using a finite volume method coupled with a hydrostatic reconstruction of the variables at the interface of each mesh cell. The
method was then validated on examples taken from the literature. Asymptotic solutions were computed to highlight the effect of
the viscous and viscoelastic source terms. Finally, the blood flow  was computed in an artery where the cross-section at rest
and the stiffness were varying. In each test case, the hydrostatic reconstruction showed good results where other simpler schemes
did not, generating spurious oscillations and
nonphysical velocities.
\end{abstract}

\section{Introduction}\label{sec:intro}

In this work we are interested in modeling and simulating blood flow in arteries with varying stiffness and cross-section. The blood flow in the main arteries of the systemic network is governed by the 3D Navier-Stokes equations which can be complicated and time-consuming to solve numerically. Fortunately, using well-known hypothesis valid in the case of blood flow in arteries (long wave approximation $D/\lambda	<<1$, axial symmetry $\partial_{\theta}=0$), this system of equations can be simplified and then integrated over the cross-section of the artery in order to obtain a 1D hyperbolic system of equations, similar to the Saint-Venant system for shallow water flows. Details on the derivation of the model are presented in section \eqref{subsec:model-mec} and can also be found in \cite{Hughes73,Olufsen00}. Finally, we are left with a set of mass and momentum conservation equations with non dimensionless variables and parameters:
 \begin{equation}
  \left\{\begin{array}{l}
          \partial_t A + \partial_x Q = 0\\
          \partial_t Q + \partial_x \left( \dfrac{Q^2}{A}+\dfrac{k}{3\sqrt{\pi}\rho}A^{3/2}\right)
          =\dfrac{A} {\sqrt{\pi}\rho} \left(\partial_x {\bf A}_0- \dfrac{2}{3}\sqrt{A}\partial_x k\right)-C_f \dfrac{Q}{A}\:,
         \end{array}\right.\label{eq:Blood-k-var}
 \end{equation}
with $A(x,t)=\pi R(x,t)^2$ the cross-section area ($R$ is the radius of the artery), $Q(x,t)=A(x,t)u(x,t)$ the discharge,
 $u(t,x)$ the mean flow velocity, $\rho$ the blood density, $C_f$ the friction coefficient and ${\bf A}_0=k\sqrt{A_0}$ with $k(x)$ the stiffness of the artery
  and $A_0(x)=\pi R_0(x)^2$ the cross-section at rest.
  
The vast majority of arteries in the systemic network are tapered, meaning that the cross-section at rest $A_0\left( x \right)$ varies throughout the length of the artery. Similarly, in the presence of arterial pathologies such as aneurysm or stenoses, the stiffness $k\left( x \right)$ of the arterial wall can vary locally. As for shallow water equations with topography, the presence of tapper or variable stiffness in an artery modifies the blood flow, and both behaviors are accounted for in \eqref{eq:Blood-k-var} through the source term ${A}  \:(\:\partial_x {\bf A}_0- 2\sqrt{A}\partial_x k/3\:)\:/{\sqrt{\pi}\rho}$. To numerically solve \eqref{eq:Blood-k-var}, it is necessary, among other things, to discretize this source term. A naive
treatment of the topography gradients will most likely generate numerical oscillations, therefore the use of the so-called well-balanced schemes is required to properly balance the fluxes and the source terms.  In the following, we will focus on a specific well-balance method, called the hydrostatic reconstruction.

We will first present the derivation of the model and its properties, then the numerical method and in particular the derivation of the well-balanced scheme applied to the case of blood flow in arteries. We will then validate our method on examples taken from
the literature and verify asymptotic behaviors of the numerical solution. Finally, we will compute the blood flow in an artery with varying cross-section and stiffness.

\section{Derivation of the 1D blood flow equations}\label{subsec:model-mec}

\begin{figure}[htb]
    \centering
    \includegraphics[width=0.6\textwidth]{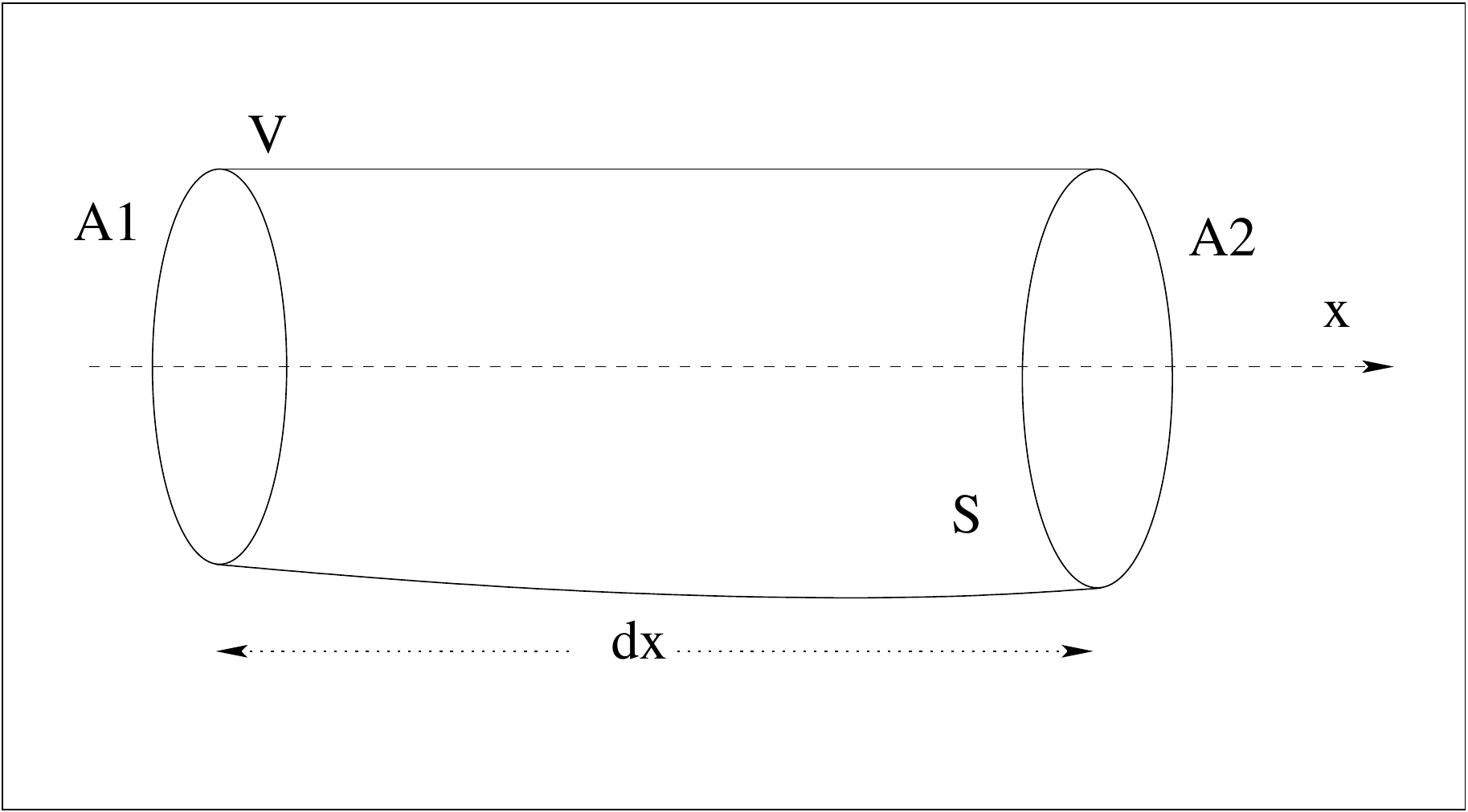}
    \caption{Control volume for integration (see text).}
    \label{fig:control}
  \end{figure}

The 1D model for blood flow equations is derived from  the conservative form of the Navier-Stokes equations for an incompressible fluid with constant
viscosity $\mu$:
\begin{eqnarray}
    \label{eq:NS1}
    \partial_t \rho + \nabla  \rho u = 0  
  \end{eqnarray}
  \begin{eqnarray}
    \label{eq:NS2}
    \partial_t \rho u + \nabla \cdot (\rho u u + p I + \tau) = 0,
  \end{eqnarray}
  where $u$ is the velocity vector, $\rho$ the density, supposed constant,
  $p$ the pressure and $\tau$ the stress tensor to be defined. Using the control volume of the Figure \ref{fig:control}, we integrate the Navier-Stokes equations over a volume $V$ of cross-section $A$ surrounded by a surface $S$ 
  ($V = S \cup  A$) and of 
  length $dz$. We define then the average velocity $U$ and pressure
  $P$ as
  \begin{eqnarray*}
    \label{eq:1}
    \{U,P\}  = \frac{1}{A} \int_{\partial A} \{u,p\}  dA.
  \end{eqnarray*}

From the mass conservation equation (\ref{eq:NS1}) we have:
\begin{eqnarray*}
    \label{eq:2}
    \int_{\partial V} (\nabla \rho u) dV = \int_{\partial S} \rho u
    \cdot n dS + \int_{\partial A}  \rho u \cdot n dA.
  \end{eqnarray*}
We then transform the volume integral using the Green (Divergence)
  theorem and writing the surface integral as $S \cup A$. The surface
  element is $dS = R d \theta dz$ and the two terms are written as
\begin{eqnarray*}
    \label{eq:3}
    \int_{\partial S} \rho u \cdot n dS = 2 \pi \int u_r|_R R dx = 2 \pi \rho \int \frac{\partial R}{\partial t} R dx =
    \rho \int \frac{\partial A}{\partial t} dx 
  \end{eqnarray*}
  and 
  \begin{eqnarray*}
    \label{eq:4}
    \int_{\partial A} \rho u \cdot n dA = (\rho A U)_1 - (\rho A
    U)_2 = \int d(\rho A U)  = \int_{\partial x} \frac{\partial \rho A U}{\partial x}, 
  \end{eqnarray*}
 We retrieve therefore the first equation of our system
  \begin{eqnarray*}
    \label{eq:5}
    \partial_t A + \partial_x (A U) = 0.
  \end{eqnarray*}
For the conservation of momentum equation \eqref{eq:NS2}, the temporal term $\partial_t \rho u$ becomes
  \begin{eqnarray*}
    \label{eq:7}
    \int_{\partial V}  \partial_t (\rho u) dV = \rho \int_{\partial V} \partial_t u dA dx = \rho \int_{\partial x} \partial_t (U A) dx 
  \end{eqnarray*}
  and the divergence term 
  \begin{eqnarray*}
    \label{eq:8}
    & \int_{\partial V} {\nabla \cdot (\rho u u + p I + \tau)} = &
    \int_{\partial S} (\rho u u + p I + \tau) \cdot n dS + \\
    &  & \int_{\partial A} (\rho u u + p I + \tau) \cdot n  dA. 
  \end{eqnarray*}
In the last two integrals the integration over the surface $S$ is
\begin{eqnarray*}
    \label{eq:9}
    \int_{\partial S} (\rho u u + p I + \tau) \cdot n  dS =  \int_{\partial S} (pn_x + \tau_{rx}) dS,
  \end{eqnarray*}
  where the term $uu dS$ tends to zero. Finally, the integration over the area $A$ gives
  \begin{eqnarray*}
    &  \int_{\partial A} (\rho u u + p I + \tau) \cdot n dA  & = [A (\rho {U^2}
    + P + \tau_{xx})]_1^2 \\
    & & = \rho \int_{\partial x} { \partial A ( {U^2} + P/\rho)
      \over \partial x } dx.
  \end{eqnarray*}
In terms of the cross-section $A$ and the flow rate $Q$, we obtained the following system of equations:
\begin{equation}
\begin{array}{l}
  \partial_t A + \partial_x Q
  = 0 \\
   \partial_t Q +  
  \partial_x {Q^2 \over A}  = -
  {A \over \rho} \partial_x P  - f_v.
 \end{array}
  \label{eq:momentum}
\end{equation}
The viscous effects are contained in $f_{v}$ which is computed by the integration of the shear stress at the wall $\tau_{rx}$ over
the internal surface $dS$. Therefore, it depends on the exact flow condition. To close the mathematical problem we need a relation between the pressure $P$ and the cross-section $A$, $P=P(A)$, called the wall or state law. For 
$f_v = C_f {Q}/{A}$ and the state law $P=P_0 + k\left( x \right)/ \sqrt{\pi} ( \sqrt{A\left( x,t \right)}-\sqrt{A_0\left( x \right)} )$, which corresponds to the elastic response of the artery, we obtain the proposed system of equations. \eqref{eq:Blood-k-var}.

\section{Conservative hyperbolic system and steady states}\label{subsec:model-prop}

Considering an artery with a constant stiffness $k$ and a variable cross-section at rest $A_0\left( x \right)$, \eqref{eq:Blood-k-var} reduces to the following system, similar to the shallow water equations with topography:
\begin{equation}
 \left\{\begin{array}{l}
         \partial_t A + \partial_x Q =0\\
         \partial_t Q + \partial_x \left(\dfrac{Q^2}{A}+\dfrac{k}{3\rho\sqrt{\pi}}A^{3/2}\right)=
         \dfrac{A}{\rho\sqrt{\pi}}\partial_x {\bf A}_0-C_f \dfrac{Q}{A}\:.
        \end{array} \right.\label{eq:Blood-k-fix}
\end{equation}
As a reminder, the shallow water system is:
\begin{equation}
 \left\{\begin{array}{l}
         \partial_t h +\partial_x q=0\\
         \partial_t q + \partial_x \left(\dfrac{q^2}{h}+\dfrac{g}{2}h^2\right)= gh \left(S_0-S_f\right)\:,
        \end{array}\right.\label{eq:SV-sys}
\end{equation}
with $h(x,t)$ the water height, $q(x,t)=h(x,t)u(x,t)$ the unit discharge, $u(x,t)$ the mean flow velocity, $g$ the constant of
gravity, $S_0=-\partial_x z$ the opposite of the slope, $z$ the topography and $S_f$ the friction term (which takes the form
 of Manning's, Stickler's, Ch\'ezy's, ... empirical friction law).
 
\subsection{Hyperbolic system}
 
  The system \eqref{eq:Blood-k-fix} can be written using the following vectorial form:
\begin{equation}
 \partial_t U + \partial_x F(U)=S(U),\label{eq:Blood-vect}
\end{equation}
where $U$ is the vector of the conservative variables, $F(U)$ is the flux:
\begin{equation}
 U=\left(\begin{array}{c}
          A\\
          Q
         \end{array}\right),\quad
F(U)=\left(\begin{array}{c}
            Q\\
            \dfrac{Q^2}{A}+\dfrac{k}{3\rho\sqrt{\pi}A^{3/2}}
           \end{array}\right)\:,\label{eq:vect-cons}
\end{equation}
and $S(U)$ is the source term, taking into account the shape of the vessel at rest $A_0\left( x \right)$ and the friction term 
\begin{equation}
 S(U)=\left(\begin{array}{c}
             0\\
             \dfrac{A}{\rho \sqrt{\pi}}\partial_x {\bf A}_0-C_f\dfrac{Q}{A}
            \end{array}\right)\:.
\end{equation}
The analogous term for the shallow water equations is the topography source term. The gradient of the flux \eqref{eq:vect-cons} can be written as the product of the Jacobian matrix $J(U)$ with the partial derivative of the vector of conservative variables $U$:
   \begin{equation}
    \partial_x F(U)=
     \left(\begin{array}{cc}
        0    							&	     1\\
        \dfrac{k\sqrt{A}}{2\rho \sqrt{\pi}}-\dfrac{Q^2}{A^2}	&	\dfrac{2Q}{A}	
     \end{array} \right).\partial_x\left(\begin{array}{c}
      A\\
      Q\end{array}\right)=J(U).\partial_x U.
   \end{equation}
 When the cross-section $A>0$, the Jacobian matrix admits two different real eigenvalues, $\lambda_1$ and $\lambda_2$:
 \begin{equation}
  \lambda_1=\dfrac{Q}{A}-\sqrt{\dfrac{k\sqrt{A}}{2\rho\sqrt{\pi}}}=u-c
  \quad \text{and} \quad
  \lambda_2=\dfrac{Q}{A}+\sqrt{\dfrac{k\sqrt{A}}{2\rho\sqrt{\pi}}}=u+c,
 \end{equation}
with $c$ the Moens-Korteweg wave propagation velocity (for the shallow water equations \eqref{eq:SV-sys}, $c=\sqrt{gh}$).
In this case, the system is said to be strictly hyperbolic, which is a generalization of the advection phenomenon
\cite{Godlewski96,Toro01,LeVeque02}: a part of the information concerning the flow propagates at the velocity $\lambda_1$ and
the other part at the velocity $\lambda_2$. For blood flow under physiological conditions, we have $\lambda_1 >0$ and $\lambda_2 <0$, hence the flow is subcritical.

\subsection{Steady states}

 Since the works of
  \cite{Bermudez94,Bermudez98} on the shallow water equations, it is well known
  that if a numerical scheme does not preserve steady states at the discrete level, spurious oscillations and artificial non zero velocities will be generated. The steady states for the system $\eqref{eq:Blood-k-fix}$ are obtained when considering a stationary flow ({\it i.e.} there is no
evolution in time) and are governed by the following equations:
\begin{equation}
 \left\{\begin{array}{l}
         \partial_x Q =0\\
         \partial_x \left(\dfrac{Q^2}{2A^2}+b\sqrt{A}-b\sqrt{A_0\left( x \right)}\right)=
         -C_f \dfrac{Q}{A^2}\:,
        \end{array} \right.\label{eq:Steady-k-fix1}
\end{equation}
 with $b=k/(\rho\sqrt{\pi})$ constant since we are considering an artery with a constant stiffness $k$. Neglecting the viscous friction effects (inviscid flow) by setting $C_f = 0$, we obtain the conservation of the discharge and Bernoulli's law for blood flow:
  \begin{equation}
   \left\{\begin{array}{l}
           Q=Q_0\\
           \dfrac{Q_0^2}{2A^2}+b\sqrt{A}-b\sqrt{A_0\left( x \right)}=cst.
          \end{array} \right..\label{eq:Steady-Bernoulli}
  \end{equation}
 In the literature \cite{Castro07,Noelle07,Thanh08,Bouchut10}, we can find well-balanced numerical methods able to preserve the following steady state:
  \begin{equation}
   \left\{\begin{array}{l}
           q=q_0\\
           \dfrac{q_0^2}{2g h^2}+h+z\left( x \right)=cst\:,
          \end{array} \right.\label{eq:Steady-Bernoulli-SW}
  \end{equation}
 which is the analogous of \eqref{eq:Steady-Bernoulli} in the case of the shallow water equations. However, these methods are complicated to handle due to the occurrence of critical points when solving \eqref{eq:Steady-Bernoulli} or \eqref{eq:Steady-Bernoulli-SW}. Therefore we chose to focus on simpler steady states that we call the rest steady states or the "man at eternal rest" equilibrium \cite{Delestre12a} by analogy with the "lake at rest" (introduced in \cite{Audusse04c}) or the hydrostatic equilibrium for the shallow water equations:
  \begin{equation}
   \left\{\begin{array}{l}
           q=u=0\\
           \partial_x (h+z\left( x \right))=\partial_x \eta = 0\:,
          \end{array}\right.\label{eq:Steady-hydrostatic-SW}
  \end{equation}
   where $\eta$ is the water level. In this case we have a hydrostatic balance between the hydrostatic pressure and the gravitational acceleration. By analogy, we have the following equilibrium for the blood flow in arteries:
  \begin{equation}
   \left\{\begin{array}{l}
           Q=u=0\\
 \partial_x \left(b\sqrt{A}-b\sqrt{A_0\left( x \right)}\right)=0\:.
          \end{array}\right.\label{eq:Steady-hydrostatic1}
  \end{equation}
 Numerical methods able to preserve at least the steady states \eqref{eq:Steady-hydrostatic1} are said to be "well-balanced" since the work of
    \cite{Greenberg96}. A wide panel of well-balanced methods has been developed for shallow water equations. Among
 others we can mention
 \cite{LeVeque98,Jin01a,Perthame01,Kurganov02,Gallouet03,Katsaounis04,Audusse04c,
 ChaconRebollo04,Noelle06,George08,Berthon12,Hou14,Bouchut04,Gosse13}.
  In \cite{Delestre12a}, we adapted the hydrostatic reconstruction introduced in \cite{Audusse04c} to the
   system with constant stiffness \eqref{eq:Blood-k-fix}.\\
   
    We will now present the hydrostatic reconstruction introduced in \cite{Audusse04c} adapted to the original system of equations \eqref{eq:Blood-k-var} with varying stiffness $k(x)$ and cross-section at rest $A_0\left( x \right)$. By a combination of the mass and momentum equations in \eqref{eq:Blood-k-var}, under some regularity assumptions, we have:
    \begin{equation}
   \partial_t u +\partial_x \left(\dfrac{u^2}{2}+\dfrac{1}{\sqrt{\pi}\rho}k\left( x \right)\sqrt{A}-\dfrac{1}{\sqrt{\pi}\rho}{\bf A}_0\left( x \right) \right)
   =-C_f \dfrac{Q}{A^2},\label{eq:k-var-combination}
  \end{equation}
with ${\bf A}_0\left( x \right)=k\left( x \right)\sqrt{A_0\left( x \right)}$. Considering a stationary flow where the viscous friction is neglected by setting $C_f=0$, we recover Bernoulli's law \eqref{eq:Steady-Bernoulli}. The notable difference is that $k$ is now a function of $x$. In the case of the "man at rest" equilibrium" (without artifacts such as \cite{Kirkman03,Munz09}) we obtain:
 \begin{equation}
  \left\{\begin{array}{l}
          Q=u=0\\
          \partial_x \left( k\left( x \right)\sqrt{A}-{\bf A}_0\left( x \right)\right) = 0\:.
         \end{array}\right.
         \label{eq:Steady-hydrostatic-kvar}
 \end{equation}
The fact that now $k$ is a function of $x$ will influence the way the well-balanced scheme is obtained. %In next section, we will first present
% the hydrostatic reconstruction applied on system \eqref{eq:Blood-k-fix} with $k=Cst$. Then we will see how it should be
%  modified to be well-balanced for system \eqref{eq:Blood-k-var} with $k\neq Cst$.
 In the following section, we will present a well-balanced scheme for system \eqref{eq:Blood-k-var}, based on the hydrostatic
  reconstruction for Saint-Venant/shallow water equations with variable pressure \cite{Bouchut04}.
  
\section{The numerical method}\label{sec:num-method}

\subsection{Numerical context}\label{subsec:num-context}

Several numerical methods have been used to solve the blood flow equations. In \cite{Stettler81}, they are solved thanks
 to the Methods of Characteristics (MOC). In \cite{Zagzoule86,Zagzoule91}, they use a conservative form of the model
 \begin{equation}
  \left\{\begin{array}{l}
          \partial_t A + \partial_x (Au)=0\\
          \partial_t u +\partial_x\left(\dfrac{u^2}{2}+\dfrac{P}{\rho} \right)=-Cf \dfrac{Q}{A^2}\:,
         \end{array}\right.\label{eq:blood-non-cons}
 \end{equation}
with the non-conserved vector $(A,u)$ and equations \eqref{eq:blood-non-cons} are solved with a
 two-step Lax-Wendroff scheme. In \cite{Stergiopulos92}, a quasi conservative form of the equations
 (with $s(U)$ a source term)
 \begin{equation}
  \left\{\begin{array}{l}
          \partial_t A + \partial_x Q=0\\
          \partial_t Q +\partial_x \left( \dfrac{Q^2}{A}\right)+\dfrac{A}{\rho}\partial_x p= s(U)\:,
         \end{array}\right.\label{eq:blood-quasi-cons}
 \end{equation}
 is solved thanks to a first order explicit in time upwind finite difference scheme. In \cite{Olufsen00}, they are the first to
  solve blood flow equations under a conservative form, thanks to a two-step Lax-Wendroff scheme.
   The solutions of the equations under the form \eqref{eq:blood-non-cons} using an upwind Discontinuous Galerkin method (used by \cite{Xiu07,Willemet11})
   and a Taylor Galerkin finite element method (also used in \cite{Martin05,Formaggia06,Melicher08}) have been compared in
   \cite{Sherwin03}. A MacCormack finite difference
    method has been applied in \cite{Fullana09} followed by \cite{Saito11}. Finite volume methods seem
    to be first used to solve these equations in \cite{Cavallini08,Cavallini10}. In \cite{Delestre12a}, a well-balanced
     finite volume method based on the hydrostatic reconstruction (introduced in \cite{Audusse04c}) is applied on system
     \eqref{eq:Blood-k-fix}, and this method is compared with a Taylor Galerkin method in \cite{Wang12}. We will present in the
     following sections the extension of the well-balanced scheme (based on an extension of the hydrostatic reconstruction)
     we have used to solve the system \eqref{eq:Blood-k-var}, which can be written under the following vectorial form
\begin{equation}
 \partial_t U + \partial_x F(U,Z)=S_1(U,Z)+S_2(U),\label{eq:blood-vec-k-var}
\end{equation}
with
\begin{equation}
U=\left(\begin{array}{c}
         A\\
         Q \end{array}\right),
         \quad Z=\left(\begin{array}{c}
        {\bf A}_0\\
         k \end{array}\right),
         \quad F(U,k)=\left(\begin{array}{c}
                  Q\\
                  \dfrac{Q^2}{A}+\dfrac{1}{3\sqrt{\pi}\rho} k A^{3/2}
                 \end{array}\right)\:,\label{eq:vec-k-var}
\end{equation}
and the source terms
\begin{equation}
 S_1(U,Z)=\left(\begin{array}{c}
                 0\\
                 \dfrac{A}{\sqrt{\pi}\rho}\left(\partial_x {\bf A}_0-\dfrac{2}{3}\sqrt{A}\partial_x k\right)
                \end{array}\right) \quad\text{and}\quad
 S_2(U)=\left(\begin{array}{c}
               0\\
               -C_f \dfrac{Q}{A}
              \end{array}\right).\label{eq:sources-k-var}
\end{equation}

\subsection{Convective step}\label{subsec:conv-step}

For the homogeneous system
\begin{equation}
\partial_t U+\partial_x F(U,Z)=0
\label{eq:blood-vec-hom}
\end{equation}
which is \eqref{eq:blood-vec-k-var} without source term, an explicit first order in time conservative scheme can be written as:
\begin{equation}
 \dfrac{U_i^{n+1}-U_i^n}{\Delta t}+\dfrac{F_{i+1/2}^n-F_{i-1/2}^n}{\Delta x}=0,\label{eq:scheme-hom}
\end{equation}
where $i$ refers to the cell $C_i=(x_{i-1/2},x_{i+1/2})=(x_{i-1/2},x_{i-1/2}+\Delta x)$ and $n$ to time $t_n$ with $t_{n+1}-t_n=\Delta t$. $U_i^n$ is an approximation of $U$:
\begin{equation*}
 U_i^n\simeq \dfrac{1}{\Delta x} \int_{x_{i-1/2}}^{x_{i+1/2}} U(x,t_n)dx\:,
\end{equation*}
and $F_{i+\frac{1}{2}}•$ is an approximation of the flux function $F(U,Z)$ at the cell interface $i+1/2$
\begin{equation*}
F_{i+1/2}^n={\bf F}(U_i^n,U_{i+1}^n,Z_i,Z_{i+1}).
\end{equation*}
This numerical flux will be detailed in subsection \ref{subsec:HLL-flux}.

\subsection{Source terms treatment}\label{subsec:source-terms}

\subsubsection{Topography source term $S_1\left( U,Z \right)$}

In the system \eqref{eq:blood-vec-k-var}, the term $S_1\left( U,Z \right)$ is involved in the steady state preservation, therefore requires a well-balanced treatment. Following a variant of the hydrostatic reconstruction \cite[p.93-94]{Bouchut04}, the variables are reconstructed locally from \eqref{eq:Steady-hydrostatic-kvar} on both sides of the interface $i+1/2$ of the cell $C_i$:
\begin{equation}
 \left\{\begin{array}{l}
         \sqrt{A_{i+1/2L}}=\max(k_i\sqrt{A_i}+\min(\Delta {{\bf A}_0}_{i+1/2},0),0)/k^*_{i+1/2}\\
	  U_{i+1/2L}=(A_{i+1/2L},A_{i+1/2L}.u_i)^t\\
	  \sqrt{A_{i+1/2R}}=\max(k_{i+1}\sqrt{A_{i+1}}-\max(\Delta {{\bf A}_0}_{i+1/2},0),0)/k^*_{i+1/2}\\
	  U_{i+1/2R}=(A_{i+1/2R},A_{i+1/2R}.u_{i+1})^t\:,
        \end{array}\right.\label{eq:rec-hydro}
\end{equation}
with $\Delta {{\bf A}_0}_{i+1/2}={{\bf A}_0}_{i+1}-{{\bf A}_0}_i=k_{i+1} \sqrt{{A_0}_{i+1}}-k_i \sqrt{{A_0}_i}$ and $k_{i+1/2}^*=\max(k_{i},k_{i+1})$.
 
In order to help the understanding of the principle of the hydrostatic reconstruction \eqref{eq:rec-hydro}, we present the hydrostatic reconstruction for the shallow water system of equations \eqref{eq:SV-sys}:
\begin{equation}
 \left\{\begin{array}{l}
         h_{i+1/2L}=\max(h_i+z_i-z_{i+1/2},0)\\
	  U_{i+1/2L}=(h_{i+1/2L},h_{i+1/2L}.u_i)^t\\
	  h_{i+1/2R}=\max(h_{i+1}+z_{i+1}-z_{i+1/2},0)\\
	  U_{i+1/2R}=(h_{i+1/2R},h_{i+1/2R}.u_{i+1})^t\:,
        \end{array}\right.\label{eq:rec-hydro-SV}
\end{equation}
with $z_{i+1/2}=\max(z_i,z_{i+1})$. The water height is reconstructed in a way that allows to have locally the hydrostatic
 equilibrium $h+z=cst$ on each side of the interface $i+1/2$. As mentioned in \cite{Audusse04c}, $\max(.,0)$ is there to ensure the positivity of the water height in case of drying and the upwind evaluation of $z_{i+1/2}$ ensures that $0\leq h_{i+1/2L}\leq h_i$ and $0\leq h_{i+1/2R}\leq h_{i+1}$, which has been proved in \cite{Audusse04c} to ensures the positivity of the water height. For blood flow equations with a constant stiffness $k$, the corresponding equilibrium writes $\sqrt{A}-\sqrt{A_0}=cst$, so $\sqrt{A}$ (respectively $-\sqrt{A_0}$) "plays the role" of $h$ (resp. $z$), thus in that case the hydrostatic reconstruction writes:
\begin{equation}
 \left\{\begin{array}{l}
         \sqrt{A_{i+1/2L}}=\max(\sqrt{A_i}-\sqrt{{A_0}_i}+\sqrt{{A_0}_{i+1/2}},0)\\
	  U_{i+1/2L}=(A_{i+1/2L},A_{i+1/2L}.u_i)^t\\
	  \sqrt{A_{i+1/2R}}=\max(\sqrt{A_{i+1}}-\sqrt{{A_0}_{i+1}}+\sqrt{{A_0}_{i+1/2}},0)\\
	  U_{i+1/2R}=(A_{i+1/2R},A_{i+1/2R}.u_{i+1})^t\:.
        \end{array}\right.\label{eq:rec-hydro-blood1}
\end{equation}
 As we have $-\sqrt{A_0}$ instead of $z$, we take $\sqrt{{A_0}_{i+1/2}}=\min(\sqrt{{A_0}_{i}},\sqrt{{A_0}_{i+1}})$, thus we have:
\begin{equation}
 \left\{\begin{array}{l}
         \sqrt{A_{i+1/2L}}=\max(\sqrt{A_i}+\min(\Delta \sqrt{A_0}_{i+1/2},0),0)\\
	  U_{i+1/2L}=(A_{i+1/2L},A_{i+1/2L}.u_i)^t\\
	  \sqrt{A_{i+1/2R}}=\max(\sqrt{A_{i+1}}-\max(\Delta \sqrt{A_0}_{i+1/2},0),0)\\
	  U_{i+1/2R}=(A_{i+1/2R},A_{i+1/2R}.u_{i+1})^t\:,
        \end{array}\right.\label{eq:rec-hydro-blood2}
\end{equation}
 with $\Delta \sqrt{A_0}_{i+1/2}=\sqrt{{A_0}_{i+1}}-\sqrt{{A_0}_i}$. We can notice that we recover reconstruction
 \eqref{eq:rec-hydro-blood2} if the stiffness $k$ is constant in reconstruction \eqref{eq:rec-hydro}. For consistency, the scheme \eqref{eq:scheme-hom} is modified as follows:
\begin{equation}
 U_i^{n+1}=U_i^n-\dfrac{\Delta t}{\Delta x}\left(F_{i+1/2L}^n-F_{i-1/2R}^n\right),\label{eq:scheme-source}
\end{equation}
where
\begin{equation*}
 \begin{array}{l}
  F_{i+1/2L}^n=F_{i+1/2}^n+S_{i+1/2L}\:,\\
  F_{i-1/2R}^n=F_{i-1/2}^n+S_{i-1/2R} \:,
 \end{array}
\end{equation*}
with
\begin{equation*}
 \begin{split}
  &F_{i+1/2}^n={\bf F}\left(U_{i+1/2L},U_{i+1/2R},k_{i+1/2}^*\right)\:,\\
  &S_{i+1/2L}=\left(\begin{array}{c}
    0\\
    {\bf P}(A_i^n,k_i)-{\bf P}(A_{i+1/2L}^n,k_{i+1/2}^*)
  \end{array}\right)\:,\\
 &S_{i-1/2R}=\left(\begin{array}{c}
    0\\
    {\bf P}(A_i^n,k_i)-{\bf P}(A_{i-1/2R}^n,k_{i-1/2}^*)
    \end{array}\right)\:,
 \end{split}
\end{equation*}
and ${\bf P}(A,k)=k\left( x \right) A^{3/2}/(3\rho \sqrt{\pi})$. Thus blood flow in a artery with varying cross-section at rest and stiffness is treated in a well-balanced way.

\subsubsection{Viscous source term $S_2\left( U \right)$}

 In system \eqref{eq:blood-vec-k-var}, the friction term $-C_f Q/A$ in $S_2\left( U \right)$ is treated
 semi-implicitly. This treatment is classical in shallow water simulations \cite{Bristeau01,Liang09b} and has proven
 efficient in blood flow simulation as well \cite{Delestre12a}. Furthermore, this treatment preserves the "dead man"
 equilibrium \eqref{eq:Steady-hydrostatic-kvar}. It consists in using first \eqref{eq:scheme-source} as a prediction step without friction, {\it i.e.}:
\begin{equation*}
 U_i^*=U_i^n-\dfrac{\Delta t}{\Delta x}\left(F_{i+1/2L}^n-F_{i-1/2R}^n\right),
\end{equation*}
then applying a semi-implicit friction correction on the predicted values ($U_i^*$):
\begin{equation*}
 A_i^*\left(\dfrac{u_i^{n+1}-u_i^*}{\Delta t}\right)=-C_f u_i^{n+1}.
\end{equation*}
Thus we get the corrected velocity $u_i^{n+1}$ and we have $A_i^{n+1}=A_i^*$.

\subsection{HLL numerical flux}\label{subsec:HLL-flux}

As presented in \cite{Delestre12a}, several numerical fluxes can be used (Rusanov, HLL, VFRoe-ncv and kinetic fluxes) for numerical simulations of blood flow in arteries. Details can be found in \cite{Bouchut04,Delestre10b,Delestre12a}. In this work we will use the HLL flux (Harten Lax and van Leer \cite{Harten83}) because it is the best compromise between accuracy and CPU time consumption (see \cite[chapter 2]{Delestre10b}). It writes:
\begin{equation*}
 {\bf F}(U_L,U_R,k^*)=\left\{\begin{array}{ll}
                               F(U_L,k^*) & \text{if}\;0\leq c_1\\
			      \dfrac{c_2 F(U_L,k^*)-c_1 F(U_R,k^*)}{c_2-c_1}
+\dfrac{c_1 c_2}{c_2-c_1}(U_R-U_L) & \text{if}\;c_1<0<c_2\\
			      F(U_R,k^*) & \text{if}\;c_2\leq 0\:,
                              \end{array}\right.
\end{equation*}
with
\[c_{1}={\inf\limits_{U=U_L,U_R}}({\inf\limits_{j\in\{1,2\}}}\lambda_{j}(U,k^*))\;
\text{and}\;c_{2}={\sup\limits_{U=U_L,U_R}}({\sup\limits_{j\in\{1,2\}}}\lambda_{j}(U,k^*)),\label{eq:flux-hll2}\]
where \(\lambda_1(U,k^*)\) and \(\lambda_2(U,k^*)\) are the eigenvalues of the system and $k^*=\max(k_L,k_R)$.

To prevent a blow up of the numerical values, we impose the following CFL (Courant, Friedrichs, Lewy) condition:
\begin{equation*}
 \Delta t\leq n_{CFL}\dfrac{\Delta x}{\max\limits_{i}(|u_i|+c_i)},
\end{equation*}
where $c_i=\sqrt{k_i \sqrt{A_i}/(2\rho\sqrt{\pi})}$ and $n_{CFL}=1$.

\section{Validation of the method}

To validate the well-balanced scheme presented in the previous sections for blood flow in arteries with varying stiffness $k\left( x \right)$
and cross-section at rest $A_0\left( x \right)$, we applied it to different test cases taken from \cite{Delestre12a}, where arteries with a varying cross-section at rest $A_0\left( x \right)$ and a constant stiffness $k$ were considered. For each of these examples, the rest equilibrium state was: $Q=0$ and $\sqrt{A}-\sqrt{A_0} = 0$
and non-reflecting boundary conditions were set at each end of the computational domain in the form of homogeneous Neumann
boundary conditions. The hydrostatic reconstruction scheme as well as a naive centered discretization of the source term were
systematically tested to clearly evaluate the benefit of using a well-balanced scheme. According to \cite{Delestre12a},
several Riemann solvers can be used, but we only display results obtained using the HLL flux. In the following, we present the numerical parameters, the analytic solution if it exists and the numerical results. For further details we refer the reader to \cite{Delestre12a}.

\subsection{"The man at eternal rest"}

We considered an artery at its equilibrium state,
where there is no flow and the radius of the cross-section at rest $R_0(x)$ varies throughout the artery, as for example in a dead man with an aneurysm. This
equilibrium state is exactly the one well-balanced methods are designed to preserve. If the topography source term is not treated correctly, non-physical velocity may be generated.\\

\begin{figure}[htp]
\begin{center}
\makebox[0.5\textwidth][c]{
\begin{minipage}[t]{.5\textwidth}
  \centering
  \includegraphics[scale=0.25,angle=0]{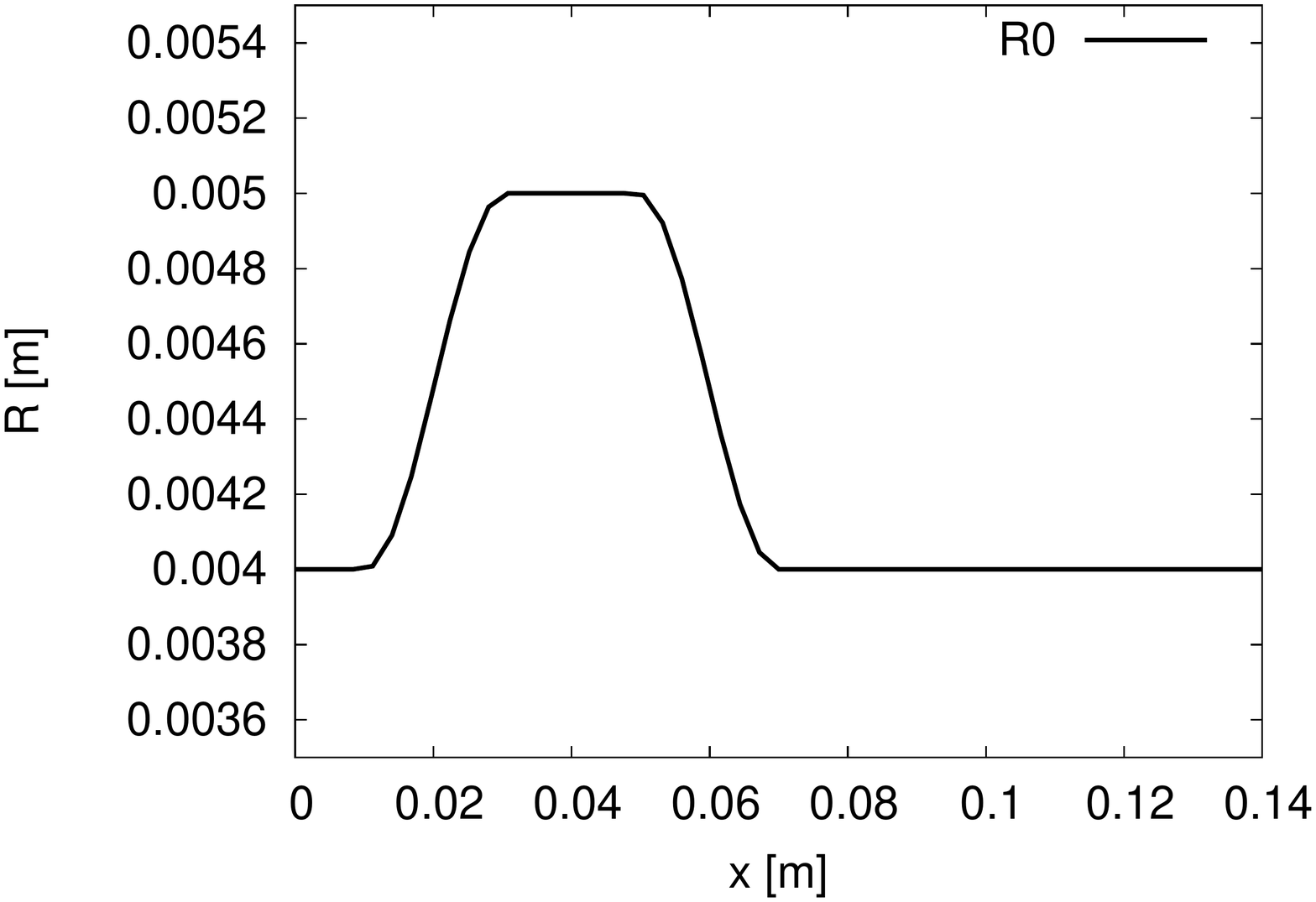}
\end{minipage}%
\begin{minipage}[t]{0.5\textwidth}
  \centering
    \includegraphics[scale=0.25,angle=0]{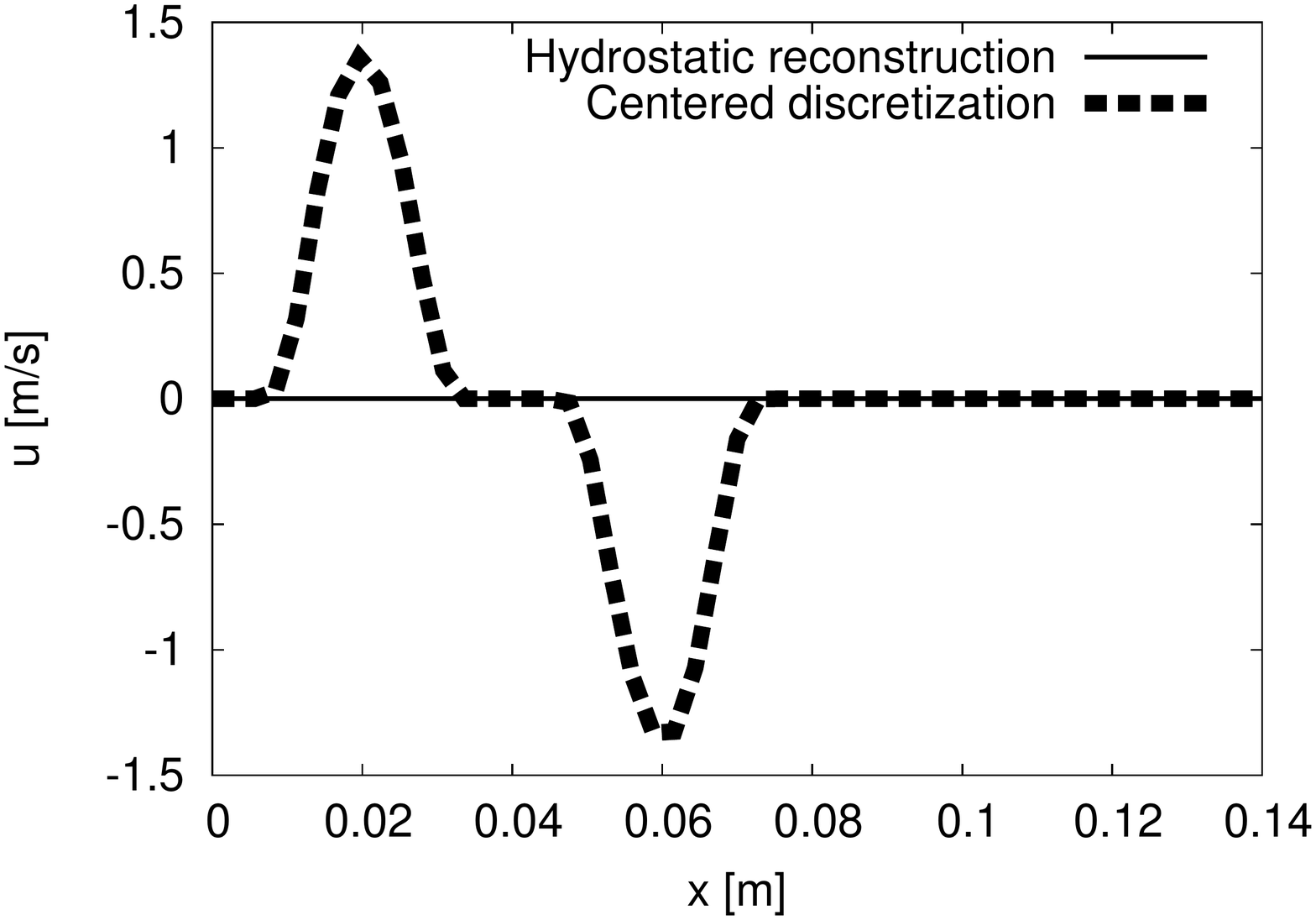}
\end{minipage}
}
\caption{The "dead man case": (Left) The radius of the artery $R_0\left( x \right)$; (Right) Comparison of the velocity at time $t=5\:s$ between an
explicit treatment of the source term (dashed line) and the hydrostatic reconstruction (full line).}
\label{fig:Dead-man-Rvar}
\end{center}
\end{figure}

We used the following numerical values:  $L = 0.14\;m$, $J = 50$ cells, $T_{end} = 5\;s$, $\rho = 1060 \;kg.m^{-3}$, $C_f = 0$ and $k = 4.0 \times 10^{8} \;Pa.m^{-1}$. We used the equilibrium state as an initial condition, setting $Q(x,0) = 0$ and:

\[
R\left( x,0 \right) = R_0\left( x \right) = 
 \left\{
 \begin{split}
 R_0 & & \text{ if } &x \in \left[0,x_1 \right] \\
 R_0 &+ \frac{\Delta R}{2}\left[ 1+\sin \left( -\frac{\pi}{2} +\pi\left( \frac{x-x_1}{x_2-x_1} \right) \right) \right] &\text{ if }& x \in \left]x_1,x_2 \right[ \\
  R_0 &+ \Delta R &\text{ if }& x \in \left[x_2,x_3 \right] \\
  R_0 &+ \frac{\Delta R}{2}\left[ 1+\cos \left( \pi\left( \frac{x-x_3}{x_4-x_3} \right) \right) \right] &\text{ if }& x \in \left]x_3,x_4 \right[ \\
  R_0 & & \text{ if } &x \in \left[x_4,L \right] \:,\\
 \end{split}
 \right.
 \] 
with $R_0 = 4.0 \times 10^{-3}\; m $, $\Delta R = 1.0 \times 10^{-3}\;m$, $x_1 = 1.0 \times 10^{-2}\;m$, $x_2 = 3.05 \times 10^{-2}\;m$,
$x_3 = 4.95 \times 10^{-2}\;m$ and $x_4 = 7.0 \times 10^{-2}\;m$. The radius at rest is plotted on Figure \ref{fig:Dead-man-Rvar} left.\\

The results obtained are presented in Figure \ref{fig:Dead-man-Rvar} right. As expected, a naive centered discretization of the topography source term results in nonphysical oscillations of the velocity $u\left( x,t \right)$, whereas the well-balanced solution
preserves the equilibrium state.\\

\subsection{The ideal "Tourniquet"}

This test case is the equivalent of the dam break problem for the Shallow Water equations (Stoker's solution in \cite{Delestre10b}). We considered an artery with a constant radius at rest $R_0$, a constant stiffness $k$ and no viscous
 friction $\left(C_f = 0\right)$, therefore the governing system of
 equations was \eqref{eq:blood-vec-hom}. Initially, a tourniquet was applied and then immediately removed. We have a Riemann problem and the method of characteristics allowed us to compute an analytic solution that we compared to the numerical solutions. This Riemann problem has been first introduced in compressible gas dynamic with the Sod tube (for further details we refer the reader to \cite{LeVeque92,Lighthill78}) and extended to blood flow in \cite{Delestre12a}.\\

We considered an artery of length $L = 8.0 \times 10^{-2}\;m$ with $x\in \left[ -\frac{L}{2},\frac{L}{2}\right]$ and used the
following numerical parameters: $J = 100$ cells, $T_{end} = 5.0 \times 10^{-3}\;s$, $\rho = 1060 \;kg.m^{-3}$ and $k = 1.0 \times 10^{7} \;Pa.m^{-1}$.  We used a perturbation of the equilibrium state as an initial condition, setting $Q(x,0) = 0$
and:

\[
A\left( x,0 \right) = 
 \left\{
 \begin{split}
 A_L =& \pi \left(R_0 + \Delta R  \right)^2 & \text{ if }& x \in \left[-\frac{L}{2},0 \right] \\
 A_R = & \pi R_0^2 &\text{ if } &x \in \left]0,\frac{L}{2} \right]\:, \\
 \end{split}
 \right.
 \] 
 with $R_0 = 4.0\times 10^{-3}\; m$ and $\Delta R = 1.0\times 10^{-3}\; m$.\\
 
\begin{figure}[htp]
\begin{center}
\makebox[0.5\textwidth][c]{
\begin{minipage}[t]{.5\textwidth}
  \centering
  \includegraphics[scale=0.25,angle=0]{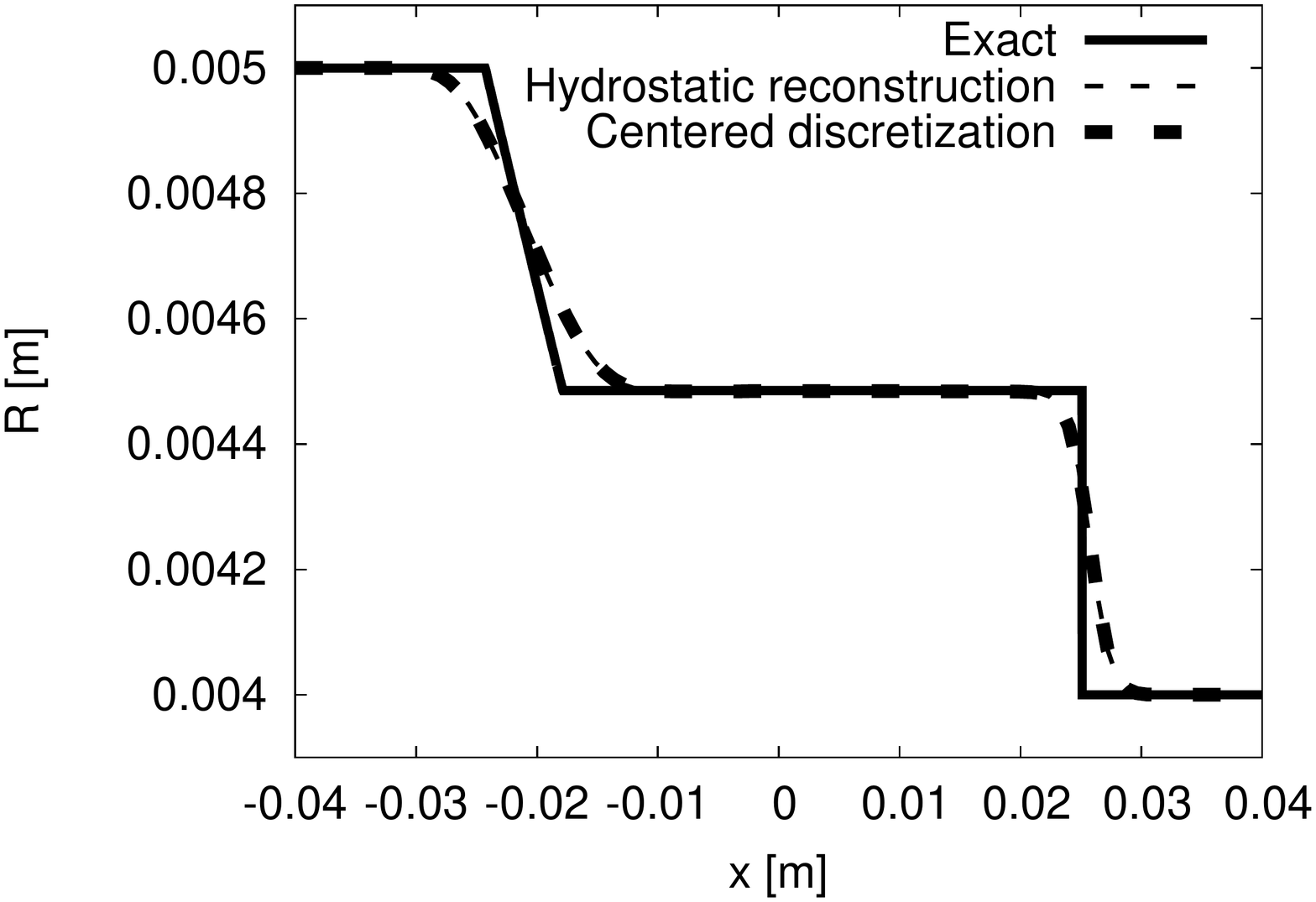}
\end{minipage}%
\begin{minipage}[t]{.5\textwidth}
  \centering
    \includegraphics[scale=0.25,angle=0]{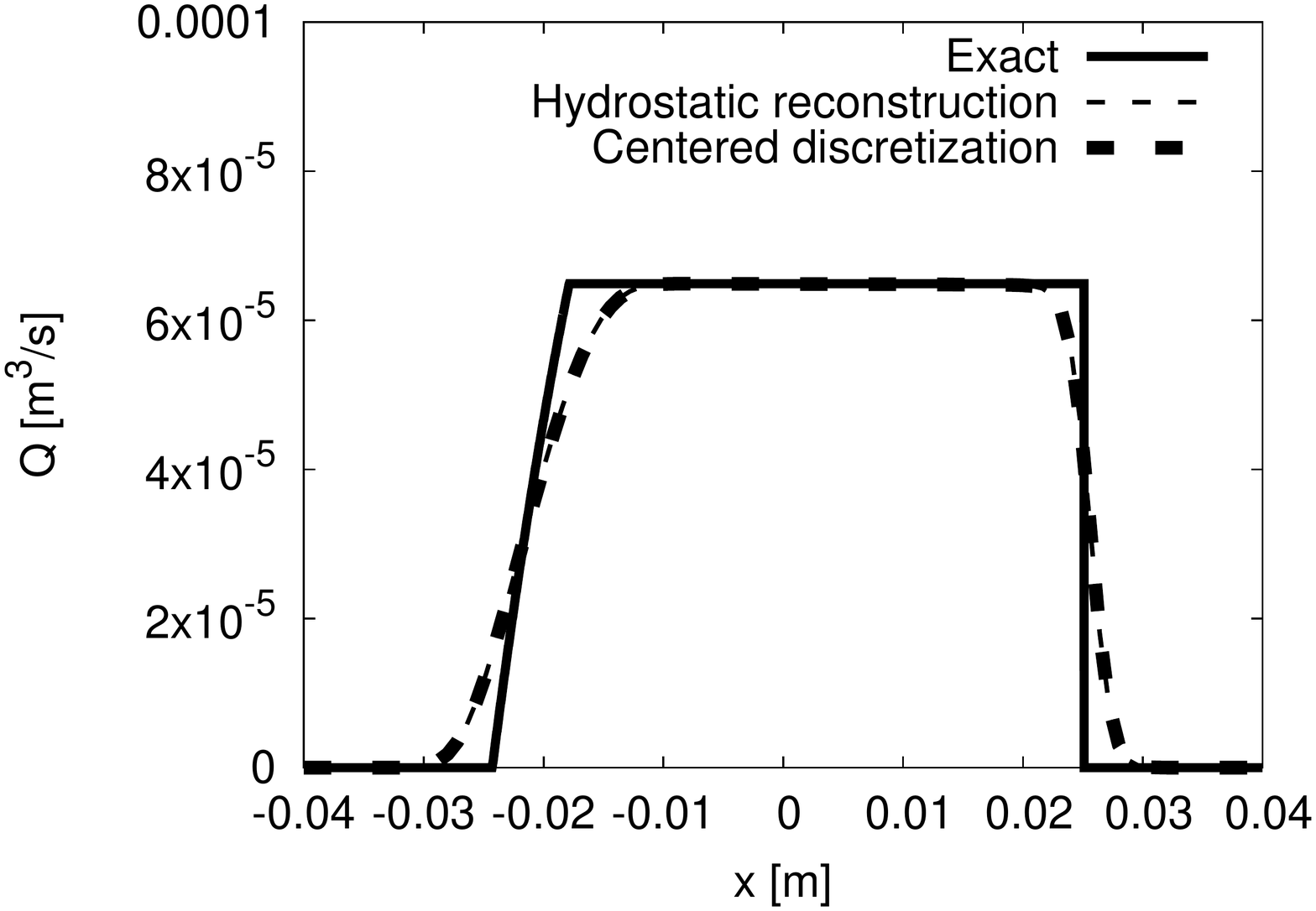}
\end{minipage}
}
\caption{The Tourniquet: (Left) Radius of the artery $R\left( x \right)$ at $t=5\times 10^{-3}\;s$; (Right) Flow rate of the artery $Q(x)$ at
$t=5\times 10^{-3}s$. Comparison between the exact analytic solution (full line) and the numerical solution obtained with an explicit
treatment of the topography source term and the hydrostatic reconstruction (dashed lines). The numerical solutions are superposed.}
\label{fig:Tourniquet}
\end{center}
\end{figure}

  The results obtained are presented in Figure \ref{fig:Tourniquet}. We can see that the numerical solution obtained with the
 well balanced scheme is in good agreement with the analytic solution presented in \cite{Delestre12a}. This is also true for
 the solution obtained using a centered discretization of the topography source term, which is superposed on the well-balanced solution, since in this case the source term is null.
  
\subsection{Wave reflection-transmission of the pulse towards a constriction}

In this section we considered the propagation of a pulse towards constriction. This configuration is an idealized representation of a transition between a parent artery and a daughter artery of smaller cross-section. We tested here
the ability of the numerical scheme to capture the propagation of a small perturbation of the equilibrium state at the beginning of an artery with a varying radius at rest $R_0(x)$. In order to accurately compute the numerical solution, the forward and backward traveling waves need to be correctly captured as well as the reflected and transmitted waves generated by the abrupt change in topography at the transition point. To test if these reflections were accurately described, we computed the analytic reflection and transmission coefficients at the transition point and compared them to the amplitude of the numerical reflected waves. For further details we refer the reader to \cite{Delestre12a}.\\

We considered an artery of length $L=0.16\; m$ and used the following numerical parameters: $J = 1500$ cells, $T_{end} = 8.0\times 10^{-3}\;s$, $\rho = 1060\;kg.m^{-3}$, $C_f = 0$ and $k = 1.0 \times 10^{8} \;Pa.m^{-1}$. The constriction was defined by
the following radius of the cross-section at rest:

\[
R_0\left( x \right) = 
 \left\{
 \begin{split}
 R_R & +\Delta R & \text{ if }& x \in \left[0,x_1 \right] \\
R_R & + \frac{\Delta R}{2} \left[1+ \cos \left( \pi\frac{x-x_1}{x_2-x_1} \right) \right] &\text{ if }& x \in \left]x_1,x_2 \right] \\
R_R & & \text{ if }& x \in \left]x_2,L \right]\:,  \\
 \end{split}
 \right.
 \] 
with $R_R = 4.0 \times 10^{-3}\; m$, $\Delta R = 1.0 \times 10^{-3}\; m$, $x_1 = \frac{19}{40}L$ and
$x_2 = \frac{L}{2}$. We set $Q(x,0)=0$ as an initial condition and we defined the initial perturbation as:

\[
R\left( x,0 \right) = 
 \left\{
 \begin{split}
 R_0(x) &\left[1+ \epsilon \sin \left( \frac{100}{20L}\pi\left(x-x_3\right) \right) \right] &\text{ if } x \in \left[x_3,x_4 \right] \\
R_0(x) & &\text{ else } \:, \\
 \end{split}
 \right.
 \] 
with $x_3 = \frac{15}{100}L < x_1$, $x_4 = \frac{35}{100}L<x_2$ and $\epsilon=5.0 \times 10^{-3}$ a small parameter ensuring
that we stayed in the range of small perturbations of the equilibrium state.\\

\begin{figure}[htp]
\begin{center}
\makebox[0.5\textwidth][c]{
\begin{minipage}[t]{.5\textwidth}
  \centering
  \includegraphics[scale=0.25,angle=0]{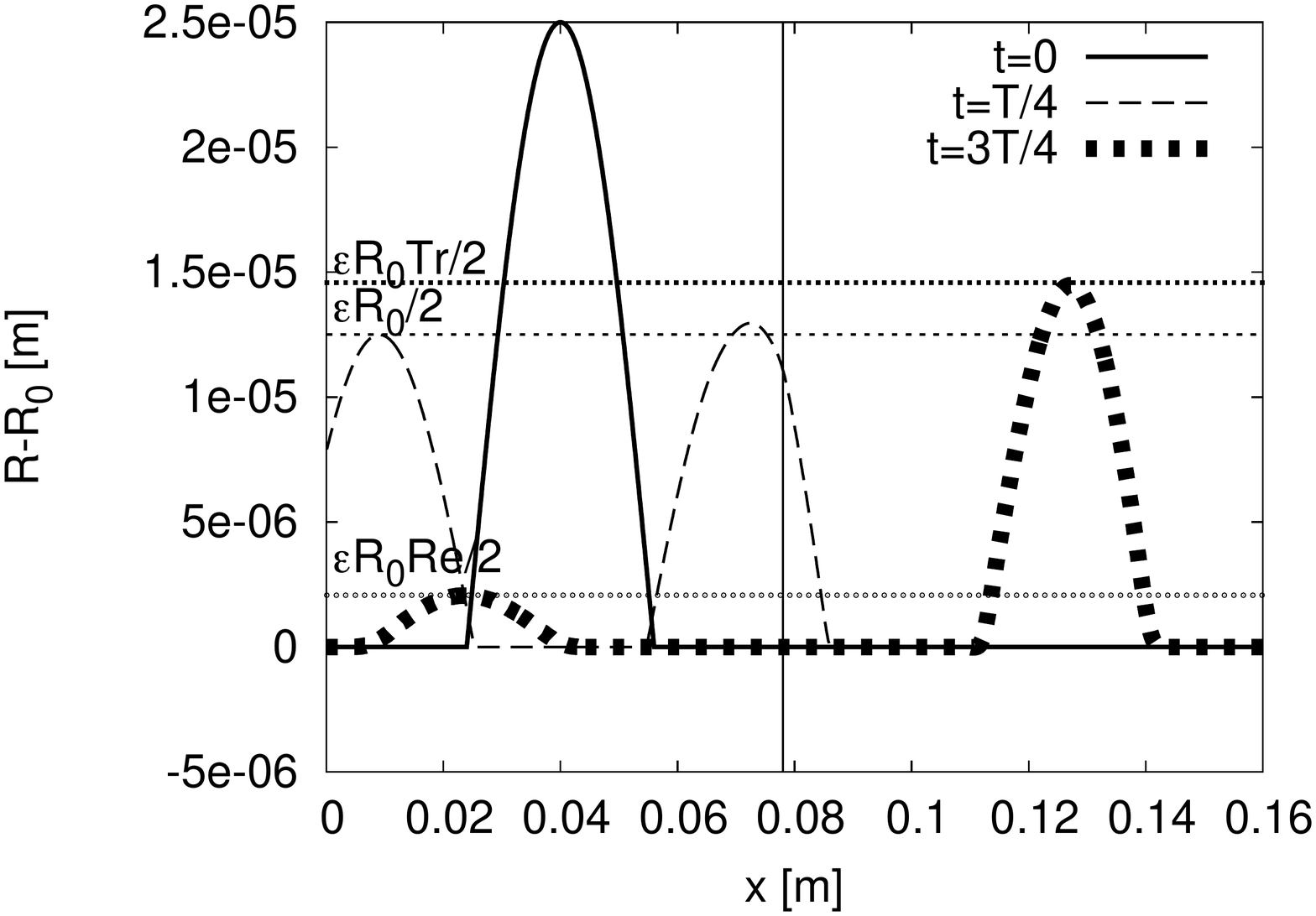}\\
\end{minipage}%
\begin{minipage}[t]{.5\textwidth}
  \centering
    \includegraphics[scale=0.25,angle=0]{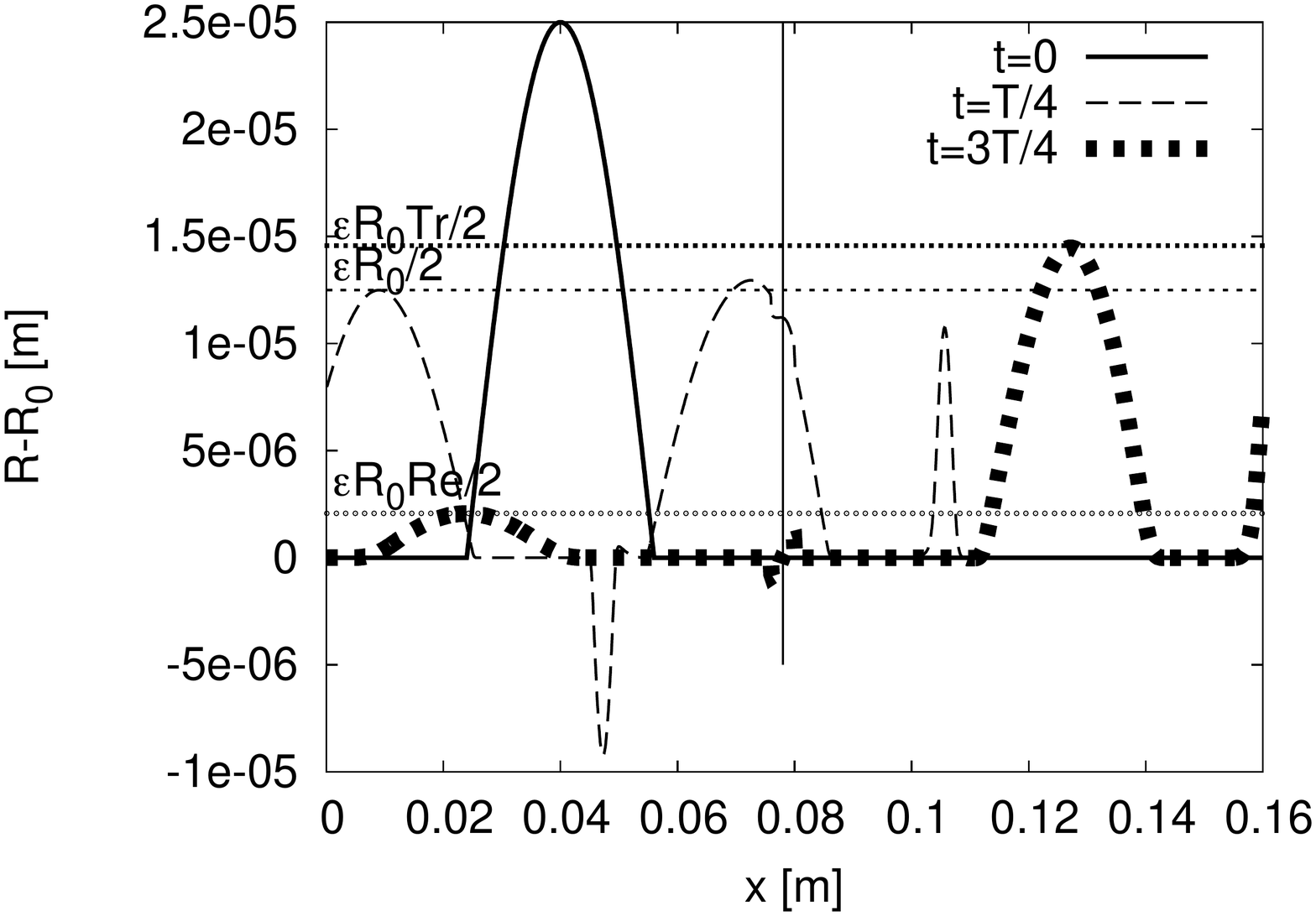}\\
\end{minipage}
}
\caption{(Left) Hydrostatic reconstruction; (Right) Centered discretization of the topography source term. $R(x)-R_0(x)$ at 3 time
steps: $t=0$, $t = \frac{T_{end}}{4}$, $t = 3\frac{T_{end}}{4}$. The straight dashed lines represent the level of the predicted reflection ($R_e$) and transmission ($T_r$) coefficients.}
\label{fig:Wave-refle-coeff}
\end{center}
\end{figure}

The numerical results are plotted in Figure \ref{fig:Wave-refle-coeff}. We can see that the propagation of the pulse as well as the wave reflections and transmissions are accurately described using the well balanced scheme (Figure \ref{fig:Wave-refle-coeff} left) whereas spurious waves appear with the centered discretization of the source term (Figure \ref{fig:Wave-refle-coeff} right).

\section{Asymptotic solutions for a uniform vessel}

In this section we studied the propagation of a pulse wave in a uniform vessel $(k=cst,A_0=cst)$ and derived
asymptotic solutions of the system of equation \eqref{eq:Blood-k-var}, following the work of Wang and al. \cite{Wang15}.
Small perturbations $\left( \epsilon \tilde{Q}, A_0 + \epsilon \tilde{A} \right)$ of the base state
$\left( Q = 0,A=A_0 \right)$ were considered, resulting in the following linearized system of equations:

\begin{equation}
\left\{
\begin{split}
          \partial_t \tilde{A}& + \partial_x \tilde{Q} =  0\\
          \partial_t \tilde{Q}& + c_0^2\partial_x \tilde{A} = -C_f \dfrac{\tilde{Q}}{A_0} \:,\\
         \end{split}\right. \label{eq:Blood-k-fix-A0-fix-linearized}
\end{equation}
where $c_0 = \sqrt{{kR_0}/\left( 2\rho \right)}$ is the Moens-Korteweg celerity.

In the following numerical examples, we only present results obtained for the hydrostatic reconstruction since
we considered a uniform vessel. The numerical parameters were defined as follows: $L =3\;m$, $R_0 =1.0\times 10^{-2} \;m$, $J = 1500$ cells, $T_{end} = 0.5\;s$, $\rho = 1060 \;kg.m^{-3}$, $\mu = 3.5\times10^{-3}\;Pa.s$ and $k = 1.0\times 10^{7}\;Pa.m^{-1}$. The parameters $C_f$ and $C_v$, respectively the viscous coefficient and the viscoelastic coefficient, were set according to the desired test case.\\

Initially, the system was at its equilibrium state $\left( Q = 0,A=A_0=\pi R_0^2 \right)$ and an inflow boundary condition was
prescribed as $Q\left(x=0,t  \right) = Q_{in}\left( t \right)$ with:

\[
\begin{split}
Q_{in}\left( t \right) = Q_c \: \text{sin}(\frac{2\pi}{T_c}t) \, H\left(-t + \frac{T_c}{2} \right) \:, \: \: t > 0\:,
\end{split}
\]
where $H(t)$ is the Heaviside function, $T_c$ the period of the sinusoidal wave and $Q_c$ the maximum amplitude of the inflow
wave. We set $Q_c =1.0 \times 10^{-6}\;m^3.s^{-1}$ and $T_c = 0.4\;s$ to insure that only small perturbations from
the equilibrium state were considered. The cross-section at the inlet $A(x=0,t)$ was reconstructed by a matching of the
outgoing characteristic, technique that takes advantage of the hyperbolic nature of the problem. A
homogeneous Neumann boundary condition was prescribed at the outlet to simplify the computation of the asymptotic solutions and
to avoid reflections.\\

\subsection{The d'Alembert equation}

Following ideas developed in \cite{Wang15}, we set $C_f = 0$ in \eqref{eq:Blood-k-fix-A0-fix-linearized} and we obtained the
d'Alembert equation, which admits the following pure wave solution $c_0\tilde{A_0}=\tilde{Q}=Q_{in}\left(x-c_0t\right)$.\\

\begin{figure}[htp]
\begin{center} \includegraphics[scale=0.25,angle=0]{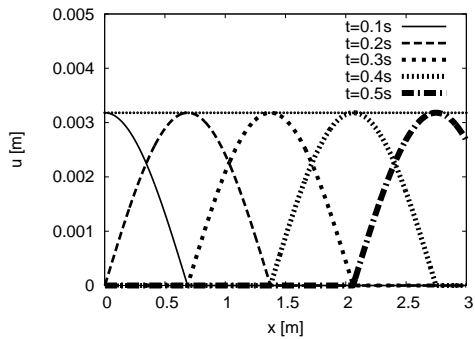}\\
\caption{Pure wave solution: $u\left( x \right)$ at time $t =\left\lbrace 0.1,0.2,0.3,0.4,0.5\right\rbrace$ for the well-balanced scheme. The straight black dotted line represents the maximum amplitude of the pure wave solution.}
\label{fig:Wave-dalembert}
\end{center}
\end{figure}

In Figure \ref{fig:Wave-dalembert}, we can see the propagation of a pulse wave without dissipation or diffusion, as predicted by the analytic solution.

\subsection{Dissipation due to the viscosity of the blood}

We also investigated the effect of the blood viscosity on the propagation of the pulse wave and set $C_f \neq 0$. Starting from the linearized system of equations \eqref{eq:Blood-k-fix-A0-fix-linearized}, we considered
the small parameter $\epsilon_f = T_c\frac{C_f}{A_0}$ and performed the change of variables $\xi = x-c_0t$ and $\tau = \epsilon_f t$ to place ourselves in the moving frame at slow times to properly capture the effects of the viscous term. The first order
solution obtained in \cite{Wang15} is:

\[
c_0 \tilde{A_0} = \tilde{Q_0} = \tilde{Q_0}\left(x-c_0t  \right)\exp{\left( -\epsilon_f \frac{t}{2T_c} \right)},
\]
where $\exp{\left( -\epsilon_f \frac{t}{2T_c} \right)}$ is the exponential envelop of the pure wave solution
$\tilde{Q_0}\left( x-c_0t \right)$. To obtain this asymptotic solution numerically, we set $C_f = 40 \pi \nu = 4.15\times 10^{-4}\;m^2.s^{-1}$, therefore $\epsilon_f = 0.53 $.\\

\begin{figure}[htp]
\begin{center} \includegraphics[scale=0.25,angle=0]{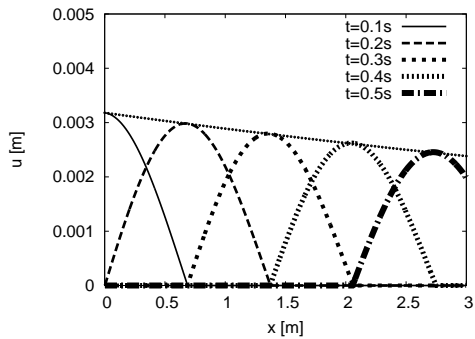}\\
\caption{Viscous damping: $u\left( x \right)$ at time $t =\left\lbrace 0.1,0.2,0.3,0.4,0.5\right\rbrace$ for the well-balanced scheme.
The straight black dotted line represents the exponential envelop of the asymptotic solution.}
\label{fig:Wave-Viscous}
\end{center}
\end{figure}

In Figure \ref{fig:Wave-Viscous}, we can see the propagation of the pulse with dissipation (or attenuation) of its amplitude due to the viscosity of the blood. The straight doted line represents the exponential envelop $\exp{\left( -\epsilon_f \frac{x}{2T_cc_0} \right)}$ computed previously and is in good agreement with the decrease in amplitude of the pulse wave. One can note that as expected, there is no diffusion, since the wavelength of the pulse does not change while it propagates in the artery.

\subsection{Diffusion due to the viscoelasticity of the arterial wall}

In this section, we set the friction coefficient to zero $(C_f=0)$ and focused on an other important characteristic of the blood flow in the arteries: the viscoelasticity of the arterial wall. We chose here to take into account this time-dependent behavior in our governing system of equations through a very simple lumped model, the Kelvin-Voigt
model, resulting in an additional parabolic term in the governing system of equations:

 \begin{equation}
  \left\{\begin{array}{l}
          \partial_t A + \partial_x Q = 0\\
          \partial_t Q + \partial_x \left( \dfrac{Q^2}{A}+\dfrac{k}{3\sqrt{\pi}\rho}A^{3/2}\right)
          =-C_f \dfrac{Q}{A} + C_v \partial_x^2 Q \:,
         \end{array}\right.\label{eq:Blood-k-fix-A0-fix-viscoelastic}
 \end{equation}
where the viscoelastic coefficient $C_\nu$ is defined as $C_{\nu} = \frac{2}{3} \frac{\phi h}{\rho R_0} = 1.57 \;m^2.s^{-1}$
with $\phi = 5000 \;Pa.s$ and $h=5.0 \times 10^{-3} \; m$.
The parabolic term was treated by performing a temporal splitting of the problem. First the purely hyperbolic problem with a non reflecting boundary condition at the outlet was
solved, and its solution was then used as an initial condition of the parabolic problem. A Crank-Nicolson scheme coupled with
homogeneous Neumann boundary conditions was than used to solve the parabolic problem.\\

To correctly capture the behavior of this new viscoelastic term, we defined a new small parameter
$\epsilon_{\nu} = \frac{C_v}{c_0^2 T_c} = 8.3 \times 10^{-2}$ and applied the same technique as in the previous section. From \cite{Wang15} we have the following first order diffusive analytic solution, which is a solution of the heat equation:

\[
\left\{
\begin{split}
&\tilde{Q_0}(\tau,\xi)  = \int_{-\infty}^{\infty} \tilde{Q_0}\left( 0,\eta \right) G\left( \tau, \xi-\eta \right)d\eta\\
&G(\tau,\xi) = \frac{1}{\sqrt{2\pi\tau c_0^2T_c}}e^{-\xi^2/\left( 2\tau c_0^2T_c \right)}\:.
\end{split}
\right.
\]

\begin{figure}[htp]
\begin{center} \includegraphics[scale=0.25,angle=0]{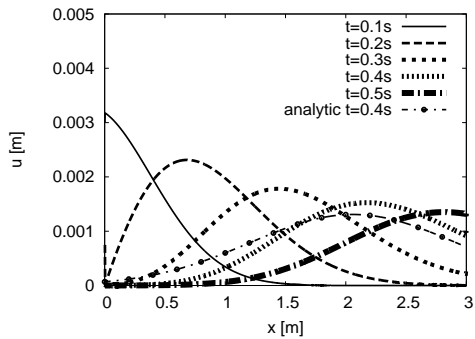}\\
\caption{Viscoelastic diffusion: $u\left( x \right)$ at time $t =\left\lbrace 0.1,0.2,0.3,0.4,0.5\right\rbrace$ for the well-balanced scheme (dashed lines). The black dotted line represents the asymptotic solution at $t=0.4 \; s$.}
\label{fig:Wave-Viscoelastic}
\end{center}
\end{figure}

The numerical results for several times and the analytic solution at $t=0.4\;s$ are presented in Figure \ref{fig:Wave-Viscoelastic}. We can see that the viscoelastic term induces a diffusion of the pulse wave, changing its wavelength, and that the numerical solution at $t=0.4\;s$ perfectly matches with the asymptotic solution at $t=0.4 \;s$.

\section{Real artery simulation}

In this section, we focused on simulating the propagation of a pulse wave in a tapered artery of length $L=3 \;m$, where the the radius of the cross section at rest $R_0(x)$ was linearly decreasing from the proximal to the distal end of the artery:

\[
R_0(x) =  
\left\{
\begin{split}
 R_L & &\text{ if }& x\in\left[0,x_1\right[\\
 R_L & -(x-x_1)\Delta R &\text{ if }& x\in\left[x_1,x_2 \right[\\
  R_L & -(x_2-x_1)\Delta R &\text{ if }& x\in\left[x_2,L \right[\:,\\
\end{split}
\right.
\]
with $R_L = 4.0 \times 10^{-3}\;m $, $\Delta R = 1.0 \times 10^{-3}\;m$, $x_1 = \frac{4}{20}L$ and $x_2=\frac{16}{20}L$. Following
\cite{Wang15}, the stiffness of the arterial wall was defined as $k(x) = \frac{4}{3}\frac{Eh}{R_0^2(x)}$ with E the Young's modulus and h the width of the arterial wall. Therefore we were in
a configuration where $R_0$ and k were varying throughout the length of the artery and if the well-balanced scheme was not used,
spurious waves might have arisen.\\

We used the following numerical parameters to mimic the geometrical and mechanical properties of a real artery:
$J = 1500$ cells, $T_{end} = 0.5\;s$, $\rho = 1060 \;kg.m^{-3}$, $\mu = 3.5\times10^{-3}\;Pa.s$, $E = 4.0 \times 10^{5} \;Pa$,
$h = 5.0 \times 10^{-4} \;m$, $C_f = 8 \pi \nu$, $\phi = 5000 \;Pa.s $ and $C_v= \frac{2}{3} \frac{\phi h}{\rho R_0}$.
We used the same initial inflow condition as for the asymptotic solutions.

\begin{figure}[htp]
\begin{center}
\makebox[0.5\textwidth][c]{
\begin{minipage}[t]{.5\textwidth}
  \centering
  \includegraphics[scale=0.25,angle=0]{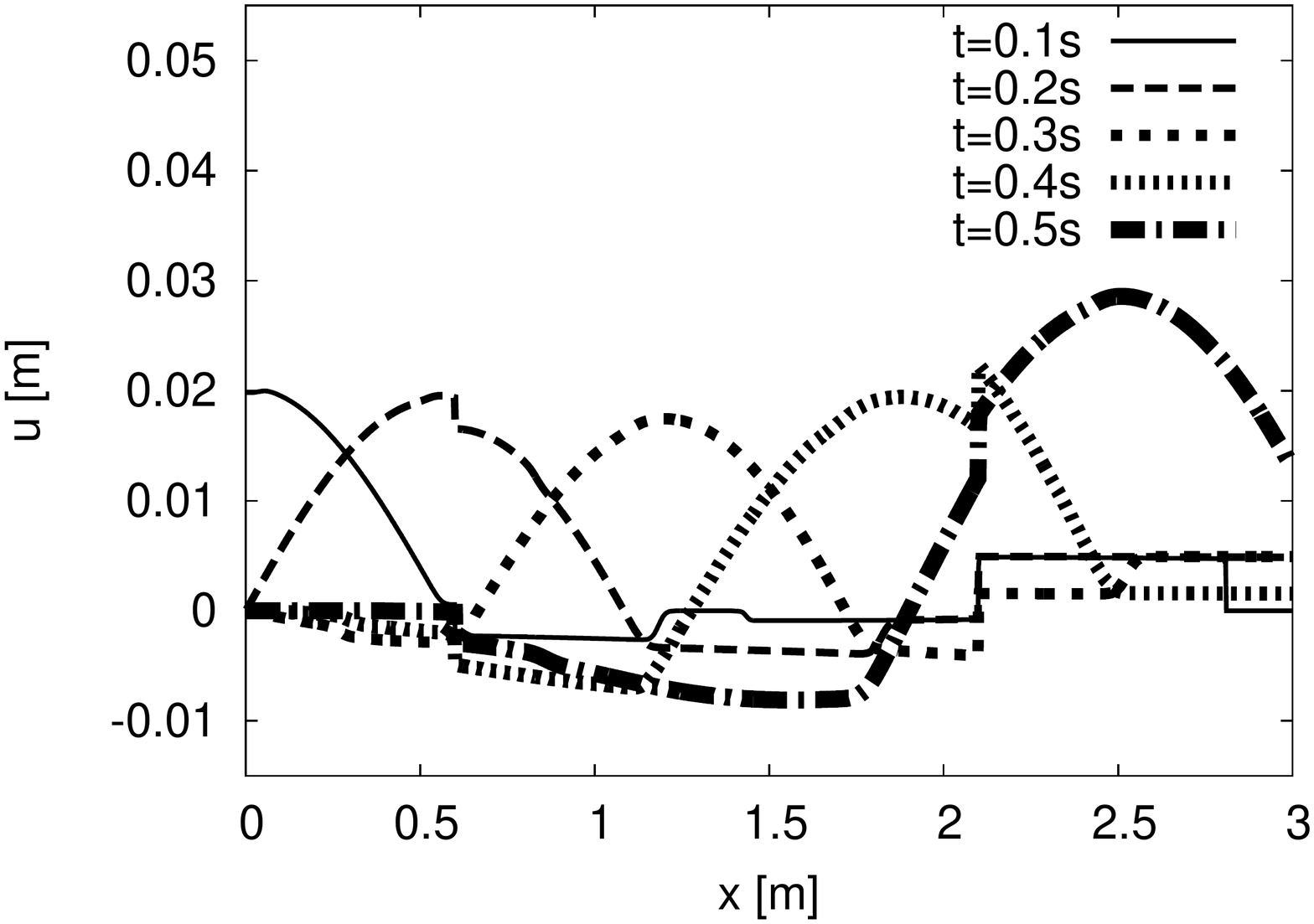}\\
\end{minipage}%
\begin{minipage}[t]{.5\textwidth}
  \centering
    \includegraphics[scale=0.25,angle=0]{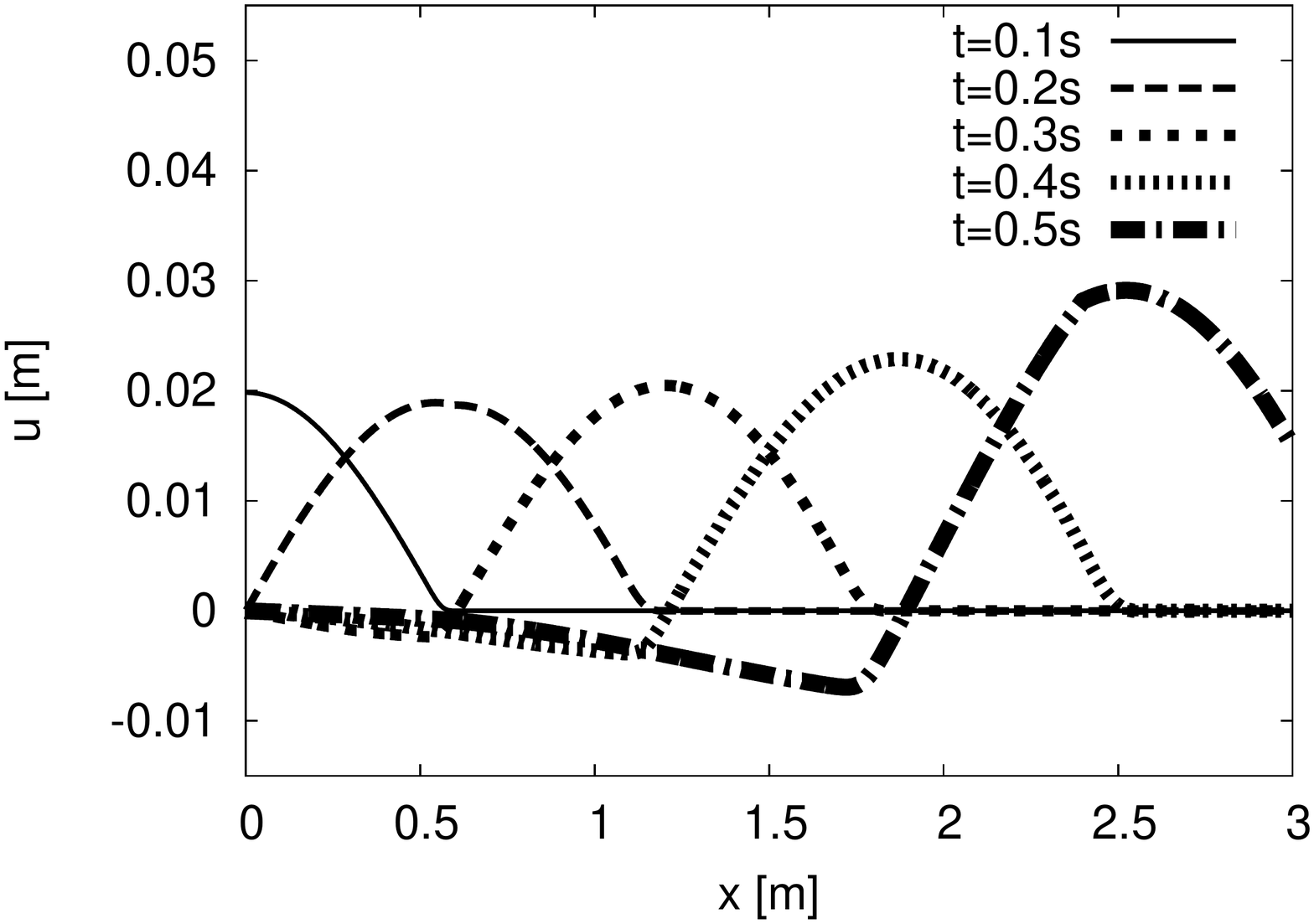}\\
\end{minipage}
}
\caption{Tapered artery - Pure wave solution: $u\left( x \right)$ at time $t =\left\lbrace 0.1,0.2,0.3,0.4,0.5\right\rbrace$ for $C_f=0$ and $C_{\nu}=0$: (Left) Centered
discretization of the topography source term; (Right) Hydrostatic reconstruction.}
\label{fig:Wave-Tapper-dAlembert}
\end{center}
\end{figure} 

\begin{figure}[htp]
\begin{center} \includegraphics[scale=0.25,angle=0]{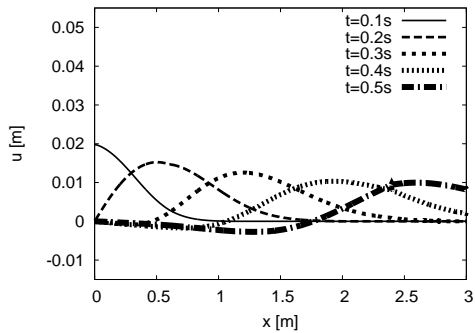}\\
\caption{Tapered artery: viscous and viscoelastic effects: $u\left( x \right)$ at time $t =\left\lbrace 0.1,0.2,0.3,0.4,0.5\right\rbrace$ for the well-balanced scheme.}
\label{fig:Wave-Tapper}
\end{center}
\end{figure} 

The results are presented in figures \ref{fig:Wave-Tapper-dAlembert} and \ref{fig:Wave-Tapper}. We can see that in the absence
of friction and viscoelastic effects (figure \ref{fig:Wave-Tapper-dAlembert}), if the well-balanced scheme is not used (figure \ref{fig:Wave-Tapper-dAlembert} left) nonphysical reflections appear. On the contrary, the well-balanced scheme provides a satisfactory numerical solution, where a continuous reflection phenomena takes place due to the tapering, resulting in a decrease of the amplitude of the backward traveling wave and an increase of the amplitude of the forward traveling wave. Indeed, in the case of a tapered artery, the transmission coefficient $T_r > 1$ and the reflection coefficient $R_e < 1$. When viscous and viscoelastic effects are taken into account (figure \ref{fig:Wave-Tapper}), all phenomena add up and we recognize the effects of the continuous reflection, the viscous dissipation
and the viscoelastic diffusion.

\section*{Conclusion and perspectives}

In this work we have presented a numerical method based on a well-balanced finite volume scheme for the blood flow equations
 with variable wall elasticity. This scheme based on an extension of the hydrostatic reconstruction gave very good results
  on several tests, for which classical methods failed. In further work, we will try to improve the accuracy of the numerical
  method by raising the order of the numerical method and to apply this method to real network modeling.
  
\section*{Acknowledgments} The first author would like to thanks the organizers of the international conference CoToCoLA
to offer the opportunity to communicate in the framework of this conference which took place in Besan\c{c}on city from
 the 9th to the 12th of February 2015.

\end{document}